\documentclass[12pt]{article}
\usepackage{amsmath}
\usepackage{amssymb}
\usepackage{latexsym}
\usepackage{amsthm}
\usepackage{tabularx}
\usepackage{booktabs}
\usepackage{mathrsfs}
\usepackage{enumitem}
\usepackage{ytableau}
\usepackage{graphicx}
\usepackage{algorithm,algorithmic}
\usepackage[colorlinks,
            linkcolor=red,
            anchorcolor=blue,
            citecolor=green]{hyperref}

\usepackage{epic}

\newtheorem{thm}{Theorem}[section]
\newtheorem{lem}[thm]{Lemma}
\newtheorem{prop}[thm]{Proposition}
\newtheorem{cor}[thm]{Corollary}

\allowdisplaybreaks

\begin{document}
\begin{center}
{\large \bf  The Smith normal form of a specialized Giambelli-type matrix}
\end{center}

\begin{center}
Alice L.L. Gao$^1$, Matthew H.Y. Xie$^2$ and Arthur L.B. Yang$^{3}$\\[6pt]

Center for Combinatorics, LPMC\\
Nankai University, Tianjin 300071, P. R. China\\[6pt]

Email: $^{1}${\tt gaolulublue@mail.nankai.edu.cn},
       $^{2}${\tt xiehongye@163.com},
       $^{3}${\tt yang@nankai.edu.cn}
\end{center}

\noindent\textbf{Abstract.}
In the study of determinant formulas for Schur functions, Hamel and Goulden introduced a class of Giambelli-type matrices with respect to outside decompositions of partition diagrams, which unify the Jacobi-Trudi matrices, the Giambelli matrices and the Lascoux-Pragacz matrices. Stanley determined the Smith normal form of a specialized Jacobi-Trudi matrix. Motivated by Stanley's work, we obtain the Smith normal form of a specialized Giambelli matrix and a specialized Lascoux-Pragacz matrix. Furthermore, we show that, for a given partition, the Smith normal form of any specialized Giambelli-type matrix can be obtained from that
of the corresponding specialization of the classical Giambelli matrix by a sequence of stabilization operations.

\noindent \emph{AMS Classification 2010:} 05E05


\section{Introduction}

Recently, there is a rising interest in the study of the Smith normal form of combinatorially defined matrices, see Stanley's survey \cite{stanley2016jcta} and references therein. The main objective of this paper is to evaluate the Smith normal form of some specializations of  Giambelli-type matrices \cite{chen2004ejc,hamel1995ejc}, which arose in the determinant formulas for Schur functions in the theory of symmetric functions. In the following, we shall explain the background and motivation of this work.

Let us first review some definitions and results on the Smith normal form.
Let $R$ be a commutative ring with identity $1$. Let $A$ be an $n\times n$ matrix over $R$.
We say that $A$ is invertible over $R$ if its determinant $\det A$ is a unit in $R$.
For an invertible $n\times n$ matrix $M$, the matrix $MA$ (or $AM$) could be obtained from $A$ by applying a sequence of elementary row (resp. column) operations, namely multiplying a row (resp. column) by a unit in $R$, or adding some multiple of a row (resp. column)
to another row (resp. column).
If $A$ can be transformed into a matrix $B$ that vanishes off the main diagonal by elementary row and column operations, then we call $B$ a diagonal form of $A$. Suppose that $A$ is of rank $r$ and the main diagonal of $B$ is $\{d_1,d_2,\ldots, d_r,0,\ldots,0\}$. If $d_k$  divides $d_{k+1}$ in $R$ for each $k:\, 1\leq k\leq r-1$, then we call $B$ a Smith normal form of $A$. Following Stanley \cite{stanley2016jcta}, we write $A\,\stackrel{\mathrm{snf}}{\rightarrow}\,(d_1,d_2,\ldots, d_n)$ to indicate that $B$ is a Smith normal form of $A$.

We would like to point out that we can study the Smith normal form in a more general setting. Indeed, $A$ is not necessarily a square matrix. For the purpose of this paper, it would be enough to consider the Smith normal form of square matrices over a principal ideal domain. It is well known that if $R$ is a principal ideal domain then $A$ always has a Smith normal form, see \cite[Theorem 17]{kuperberg2002ejc}, and moreover, the diagonal entries of the Smith normal form are unique up to multiplication by units.
The following result provides a useful formula for the Smith normal form over a principal ideal domain, see \cite[Theorem 2.4]{stanley2016jcta}.

\begin{thm}\label{thm-minor}
Suppose that $R$ is a principal ideal domain and $A$ is an $n\times n$ matrix over $R$.
If $A\,\stackrel{snf}{\rightarrow}\,(d_1,d_2,\ldots, d_n)$, then, for each $k:\, 1\leq k\leq n$, the product $d_1d_2\cdots d_k$ is equal to the greatest common divisor of all $k\times k$ minors of $A$. (By convention, we set the greatest common divisor to be $0$ if all $k\times k$ minors are $0$.)
\end{thm}

This paper is motivated by Stanley's work on the Smith normal form of some specializations of the Jacobi-Trudi matrix \cite{stanley2015arxiv,stanley2016jcta}.
To state Stanley's result, we shall first recall some definitions on partitions.  Let $n$ be a nonnegative integer. By a {partition} of $n$ we mean a tuple $\lambda=(\lambda_1,\lambda_2,\ldots,\lambda_k)$ of nonnegative integers such that $\lambda_1\geq\lambda_2\geq\ldots\geq\lambda_k \geq 0$ and $\sum_{i=1}^k\lambda_i=n$. The nonzero entries $\lambda_i$ are called the {parts} of $\lambda$, and the number of parts is called the {length} of $\lambda$, denoted by $\ell(\lambda)$.
Each partition is associated to a left justified array of cells, called the Ferrers or Young diagram of $\lambda$.
Given two partitions $\lambda,\mu$ such that $\lambda_i\geq\mu_i$ for all $i\geq 1$ (denoted by $\mu \subseteq \lambda$), let $\lambda/\mu$ denote the skew partition corresponding to the diagram of $\lambda$ with the diagram of $\mu$ removed from its upper left-hand corner.
Given a partition $\lambda=(\lambda_1,\lambda_2,\ldots,\lambda_{\ell(\lambda)})$,  the corresponding diagram has $\lambda_i$ cells in the $i$-th row. Here we number the rows from top to bottom and the columns from left to
right. The cell in the $i$-th row and $j$-th column is denoted by $(i, j)$.
The {content} of $(i,j)$ is defined to be $j-i$, denoted by $c(i,j)$. The hook-length of $(i,j)$, denoted by $h(i,j)$, is defined
to be the number of cells directly to the right or directly below $(i,j)$, counting $(i,j)$ itself once.
We say that the rank of $\lambda$ is $r$, denoted $\mathrm{rank}(\lambda)$, if $(r,r)\in \lambda$ but $(r+1,r+1)\not\in \lambda$. For each $1\leq i\leq \mathrm{rank}(\lambda)$, let $D_i(\lambda)$ denote the $i$-th diagonal hook of $\lambda$. For the convenience, we set $D_i(\lambda)=\emptyset$
if $i>\mathrm{rank}(\lambda)$.

Given a partition $\lambda$, the Jacobi-Trudi matrix $\mathrm{JT}_{\lambda}$ is defined by
\begin{align*}
\mathrm{JT}_{\lambda}=(h_{\lambda_i-i+j})_{i,j=1}^{\ell(\lambda)},
\end{align*}
where $h_i$ is the $i$-th complete symmetric function in the variables $x_1,x_2,\ldots,$ with the convention that $h_0=1$ and $h_i=0$ for $i<0$. Let $s_{\lambda}$ denote the Schur function indexed by $\lambda$. The well known Jacobi-Trudi identity states that
\begin{align*}
s_{\lambda}=\det \mathrm{JT}_{\lambda}.
\end{align*}

For a symmetric function $f$, let $\varphi_tf$ denote the specialization $f(1^t)$, that is, set $x_1=\cdots =x_t=1$ and all other $x_i=0$ in $f$.
The hook-content formula \cite[Corollary 7.21.4]{stanley1999ec2} tells us that
\begin{align}
\varphi_t\,s_{\lambda}=\prod_{(i,j)\in \lambda}\frac{(t+c(i,j))}{h(i,j)}. \label{eq-hook-content}
\end{align}
In particular,
\begin{align*}
\varphi_t\,h_i=\binom{t+i-1}{i},
\end{align*}
which is a polynomial in $t$ of degree $i$ with rational coefficients. Since $\mathbb{Q}[t]$ (the ring of polynomials in $t$ with rational coefficients) is a principal ideal domain, and $\varphi_t\,s_{\lambda}$ factors a lot over $\mathbb{Q}[t]$, Stanley was motivated to study the Smith normal form of
the specialized Jacobi-Trudi matrix
\begin{align*}
\varphi_t\,\mathrm{JT}_{\lambda}=\left(\varphi_t\,h_{\lambda_i-i+j}\right)_{i,j=1}^{\ell(\lambda)},
\end{align*}
and obtained the following result.

\begin{thm}[{\cite[Theorem 5.3]{stanley2016jcta}}]\label{thm-snf-jt}
Suppose that $\varphi_t\,\mathrm{JT}_{\lambda}\,\stackrel{\mathrm{snf}}{\rightarrow}\, (d_1,d_2,\ldots,d_{\ell(\lambda)})$ over $\mathbb{Q}[t]$.
Then, for $1\leq k\leq \ell(\lambda)$, we can take
\begin{align*}
d_k=\prod_{(i,j)\in D_{\ell(\lambda)-k+1}}(t+c(i,j)).
\end{align*}
\end{thm}

There are also other determinant formulas for the Schur function $s_{\lambda}$. It is natural to study the Smith normal form of such matrices.
The first candidate to come to mind is the Giambelli matrix, since, as we see above, the diagonal hooks play an important role in the Smith normal form of a specialized Jacobi-Trudi matrix. The Giambelli matrix can be easily described by using the Frobenius notation of partitions.
Suppose that $\lambda$ is of rank $r$. For $1\leq i\leq r$, let $\alpha_i$ denote the number of cells
directly to the right of $(i,i)$, and let $\beta_i$ denote the number of cells
directly below $(i,i)$. Then we can denote that partition $\lambda$ by
$(\alpha|\beta)=(\alpha_1,\ldots,\alpha_r|\beta_1,\ldots,\beta_r)$, called the Frobenius notation of $\lambda$.
A moment's thought shows that $\alpha_1>\alpha_2\cdots>\alpha_r\geq 0$ and $\beta_1>\beta_2\cdots>\beta_r\geq 0$.
Recall that the {Giambelli matrix} $\mathrm{G}_{\lambda}$ is defined by
\begin{align*}
\mathrm{G}_{\lambda}=(s_{(\alpha_i | \beta_j)})_{i,j=1}^{r}.
\end{align*}
The Giambelli identity asserts that $s_{\lambda}=\det\,\mathrm{G}_{\lambda}$, see \cite{giambelli1903at}.
Now consider the specialization $\varphi_t\,\mathrm{G}_{\lambda}$.
Following Stanley's proof of Theorem \ref{thm-snf-jt}, we obtain the following result without much difficulty by using Theorem \ref{thm-minor}.

\begin{thm}\label{thm-snf-g}
Suppose that $\varphi_t\,\mathrm{G}_{\lambda}\,\stackrel{\mathrm{snf}}{\rightarrow}\, (d_1,d_2,\ldots,d_{r})$  over $\mathbb{Q}[t]$, where $r=\mathrm{rank}(\lambda)$.
Then, for $1\leq k\leq r$, we can take
\begin{align*}
d_k=\prod_{(i,j)\in D_{r-k+1}}(t+c(i,j)).
\end{align*}
\end{thm}

As we see above, the Giambelli identity gives a determinantal expression for $s_{\lambda}$ involving hook functions.
Lascoux and Pragacz \cite{lascoux1988ejc} showed that such a determinant $\det\,\mathrm{G}_{\lambda}$ can be transformed into a determinant of ribbon functions. Recall that a ribbon (or a border strip) is a connected skew diagram with no $2\times 2$ block of cells.
We say that a skew diagram is {connected} if it can be regarded as a union of an edgewise connected set of cells, where two cells are said to be {edgewise connected} if they share a common edge.
The rim of a diagram is the maximal outer ribbon of the diagram. Given a partition $\lambda$ with rank $r$, we can peel its diagram off into successive rims $\theta_1,\theta_2,\ldots,\theta_r$ beginning from the outside. It is clear that each $\theta_i$ is cut by the diagonal into three disjoint parts: $\theta_i^+$, $\Box_i$ and $\theta_i^-$, which are respectively the cells of $\theta_i$ strictly above the diagonal, the diagonal cell, and the cells strictly below the diagonal. Given two ribbons $\theta_i$ and $\theta_j$, let $\theta_i^+ \& \theta_j^-$ denote the ribbon obtained by replacing the lower part $\theta_i^-$ in $\theta_i$ by $\theta_j^-$.
Let
\begin{align*}
\mathrm{LP}_{\lambda}=(s_{\theta_i^+ \& \theta_j^-})_{i,j=1}^{r},
\end{align*}
where $s_{\theta_i^+ \& \theta_j^-}$ denotes the skew Schur function corresponding to the ribbon $\theta_i^+ \& \theta_j^-$.
Lascoux and Pragacz \cite{lascoux1988ejc} proved that $s_{\lambda}=\det\,\mathrm{LP}_{\lambda}$.
For this reason, we call $\mathrm{LP}_{\lambda}$ the Lascoux-Pragacz matrix with respect to $\lambda$.
Now consider the specialization $\varphi_t\,\mathrm{LP}_{\lambda}$ and we get its Smith normal form as follows.

\begin{thm}\label{thm-snf-lp}
Suppose that $\varphi_t\,\mathrm{LP}_{\lambda}\,\stackrel{\mathrm{snf}}{\rightarrow}\, (d_1,d_2,\ldots,d_{r})$  over $\mathbb{Q}[t]$, where $r=\mathrm{rank}(\lambda)$.
Then, for $1\leq k\leq r$, we can take
\begin{align*}
d_k=\prod_{(i,j)\in D_{r-k+1}}(t+c(i,j)).
\end{align*}
\end{thm}

Looking at the diagonal entries of the Smith normal forms in Theorems \ref{thm-snf-jt}, \ref{thm-snf-g} and \ref{thm-snf-lp}, we see
some coincidence on non-trivial entries. Since the Jacobi-Trudi matrix, the Giambelli matrix and the Lascoux-Pragacz matrix can be considered as special cases of Giambelli-type matrices introduced by Hamel and Goulden \cite{hamel1995ejc}, we were inspired to study
the Smith normal form of any specialized Giambelli-type matrix. These Giambelli-type matrices are associated with outside decompositions of the partition diagrams. Recall that, for a given partition $\lambda$, an outside decomposition $\Pi$ of $\lambda$ is a partition of the cells of $\lambda$ into pairwise disjoint border strips such that the starting cell of each strip in the decomposition is on the left or bottom perimeter and the ending cell is on the right or top perimeter. Then for each outside decomposition $\Pi$ of $\lambda$, Hamel and Goulden \cite{hamel1995ejc} defined a matrix $\mathrm{M}_{\lambda}(\Pi)$ and showed that $s_{\lambda}=\det\,\mathrm{M}_{\lambda}(\Pi)$. The matrix $\mathrm{M}_{\lambda}(\Pi)$ is called a Giambelli-type matrix. The construction of $\mathrm{M}_{\lambda}(\Pi)$ will be illustrated in Section \ref{sec-gt}.
We next consider the specialization of $\varphi_t\,\mathrm{M}_{\lambda}(\Pi)$, and obtain the following result.

\begin{thm}\label{thm-snf-gg}
Suppose that $\mathrm{M}_{\lambda}(\Pi)$ is of order $m$ and $\varphi_t\,\mathrm{M}_{\lambda}(\Pi)\,\stackrel{\mathrm{snf}}{\rightarrow}\, (d_1,d_2,\ldots,d_{m})$  over $\mathbb{Q}[t]$.
Then, for $1\leq k\leq m$, we can take
\begin{align*}
d_k=\prod_{(i,j)\in D_{m-k+1}}(t+c(i,j)).
\end{align*}
\end{thm}

It is known that if all strips in $\Pi$ are horizontal, then $\mathrm{M}_{\lambda}(\Pi)$ is $\mathrm{JT}_{\lambda}$; if all strips in $\Pi$ are hooks, then $\mathrm{M}_{\lambda}(\Pi)$ is just $\mathrm{G}_{\lambda}$; and if $\Pi$ is the rim ribbon decomposition of $\lambda$, then $\mathrm{M}_{\lambda}(\Pi)$ becomes $\mathrm{LP}_{\lambda}$.

The rest of the paper is organized as follows. To be self-contained, we will provide a detailed review of Stanley's proof of Theorem \ref{thm-snf-jt} in Section \ref{sec-jt}. Following Stanley's proof,  we will give a proof of Theorem \ref{thm-snf-g} in Section \ref{sec-g}, and give a proof of Theorem \ref{thm-snf-lp} in Section \ref{sec-lp}.
As will be shown later, we need not use the Littlewood-Richardson rule for the evaluation of the Smith normal form of a specialized Giambelli matrix.
However, for a general outside decomposition $\Pi$ of $\lambda$, it is difficult to give a proof of Theorem \ref{thm-snf-gg} along the lines of Stanley's proof of Theorem \ref{thm-snf-jt}.
In Section \ref{sec-gt} we will instead prove Theorem \ref{thm-snf-gg} based on the stable equivalence of Giambelli-type matrices, which was established by Chen and Yang \cite{chen2004ejc} in answer to a question of Kuperberg \cite{kuperberg2002ejc}. Lastly, in Section \ref{sec-gq}, we will follow Stanley to give two kinds of $q$-analogues of the Smith normal form of the specialized Giambelli-type matrices.

\section{The Jacobi-Trudi matrix}\label{sec-jt}

The aim of this section is to give an overview of Stanley's proof of Theorem \ref{thm-snf-jt}.

In order to use Theorem \ref{thm-minor} to prove Theorem \ref{thm-snf-jt}, we need to consider $k\times k$ minors of
$\varphi_t\,\mathrm{JT}_{\lambda}$ for $1\leq k\leq \ell(\lambda)$.
Each minor is either zero or a skew Schur function (under the specialization $\varphi_t$) for some skew partition $\rho/\sigma$ by the following Jacobi-Trudi identity for skew Schur functions, see \cite[Theorem 7.16.1]{stanley1999ec2}:
\begin{align*}
s_{\rho/\sigma}=\det(h_{\rho_i-\sigma_j-i+j})_{i,j=1}^{\ell(\rho)}.
\end{align*}

For the square submatrix $N_k$ of $\varphi_t\,\mathrm{JT}_{\lambda}$ with row indices $1\leq i_1<i_2<\cdots<i_k\leq \ell(\lambda)$ and column indices $1\leq j_1<j_2<\cdots<j_k\leq \ell(\lambda)$, we may take
\begin{align}
\rho_1&=\Xi+\lambda_{i_1}-i_1-k+1,\,\rho_2=\Xi+\lambda_{i_2}-i_2-k+2,\,\ldots,\,
\rho_k=\Xi+\lambda_{i_k}-i_k,\nonumber\\[5pt]
\sigma_1&=\Xi-j_1-k+1,\,
\sigma_2=\Xi-j_2-k+2,\,\ldots,\,\sigma_k=\Xi-j_k,\label{eq-part}
\end{align}
where $\Xi$ can be any positive integer such that all parts of $\rho$ and $\sigma$ are nonnegative.
If $\sigma \not\subseteq \rho$, namely $\sigma_l>\rho_l$ for some $1\leq l\leq k$, we set $s_{\rho/\sigma}=0$ for the convenience though $\rho/\sigma$ is not a valid skew partition. Now we have $\varphi_t\,s_{\rho/\sigma}=\det N_k$ in any case.

Let $M_k$ denote the square submatrix consisting of the last $k$ rows and first $k$ columns of $\varphi_t\,\mathrm{JT}_{\lambda}$.
By \eqref{eq-part}, it is routine to verify that $M_k$ is the specialized Jacobi-Trudi matrix $\varphi_t\,\mathrm{JT}_{\mu^{(k)}}$, where
\begin{align}\label{eq-part-mk}
\mu^{(k)}=(\lambda_{\ell(\lambda)-k+1}-\ell(\lambda)+k,\lambda_{\ell(\lambda)-k+2}-\ell(\lambda)+k,\ldots,\lambda_{\ell(\lambda)}-\ell(\lambda)+k).
\end{align}
Note that some components of $\mu^{(k)}$ might be negative. Let $\mu^{(k,+)}$ denote the partition consisting of nonnegative parts of $\mu^{(k)}$.
If all parts of $\mu^{(k)}$ are negative, define $\mu^{(k,+)}$ to be the empty partition.  With this convention, it is readily to see that $\mu^{(k,+)}$ admits the following hook decomposition
\begin{align}\label{eq-decomposition}
\mu^{(k,+)}=\bigcup_{l=1}^kD_{\ell(\lambda)-l+1}.
\end{align}

Let $\langle\,,\,\rangle$ denote the scalar product on the ring of symmetric functions by requiring that the Schur functions form an orthonormal basis. Stanley \cite{stanley2015arxiv} noted the following property.

\begin{lem}\label{lem-stanley} Suppose that $\det N_k\neq 0$, and $\rho/\sigma, \mu^{(k,+)}$ are defined as above.
Then there is a subdiagram $\nu$ (of an ordinary partition) of $\rho/\sigma$ containing $\mu^{(k,+)}$, and all other cells of
$\rho/\sigma$ are to the left of $\nu$. Furthermore, we have $\mu^{(k,+)}\subseteq \tau$ if $\langle s_{\rho/\sigma},\, s_{\tau} \rangle\neq 0$.
\end{lem}

\proof
To prove the first part, it suffices to show that, for any $1\leq l\leq \ell(\mu^{(k,+)})$,
 $$\rho_l-\sigma_1\geq \mu^{(k,+)}_l=\mu^{(k)}_l,$$
 that is,
 $$(\Xi+\lambda_{i_l}-i_l-k+l)-(\Xi-j_1-k+1)\geq \lambda_{\ell(\lambda)-k+l}-\ell(\lambda)+k.$$
  Equivalently, we only need to prove that
 $$
 \lambda_{i_l}-i_l+l+(j_1-1)\geq \lambda_{\ell(\lambda)-k+l}-\ell(\lambda)+k.
 $$
 Since $j_1\geq 1$ and $i_l\leq \ell(\lambda)-k+l$ (as $M_k$ consists of the last $k$ rows and first $k$ columns of $\varphi_t\,\mathrm{JT}_{\lambda}$), we obtain the desired result.

 For the second part of the lemma, the statement that $\langle s_{\rho/\sigma},\, s_{\tau} \rangle\neq 0$ is equivalent to $c_{\sigma\tau}^{\rho}\neq 0$, where
 $c_{\sigma\tau}^{\rho}$ is a Littlewood-Richardson coefficient \cite[(7.64)]{stanley1999ec2}. The celebrated Littlewood-Richardson rule states that $c_{\sigma\tau}^{\rho}$
 is equal to the number of semistandard Young tableaux of shape $\rho/\sigma$ and content $\tau$ whose reverse reading words are lattice permutations.
 By the first part of the lemma,
 such a tableau must have the last $\mu^{(k,+)}_l$ entries in row $l$ equal to $l$ for any $1\leq l\leq \ell(\mu^{(k,+)})$. Hence $\mu^{(k,+)}_l\leq \tau_l$ for any $l$, namely $\mu^{(k,+)}\subseteq \tau$.
 \qed

Stanley \cite{stanley2015arxiv} further obtained the following results.

\begin{prop}\label{prop-gcd}
\begin{itemize}
\item[(C1)] The matrix $\varphi_t\,\mathrm{JT}_{\lambda}$ has a $k\times k$ submatrix with determinant equal to
$\varphi_t\,s_{\mu^{(k,+)}}$.

\item[(C2)] Every $k\times k$ minor of $\varphi_t\,\mathrm{JT}_{\lambda}$ is divisible by $\varphi_t\,s_{\mu^{(k,+)}}$
in the ring $\mathbb{Q}[t]$.
\end{itemize}
\end{prop}

\proof First we prove (C1).
Suppose that $\det M_k\neq 0$, in which case $\mu^{(k,+)}=\mu^{(k)}$. Now $M_k$ serves as a candidate since $\det M_k=\varphi_t\,s_{\mu^{(k)}}$.
If $\det M_k=0$, then some components of $\mu^{(k)}$ in \eqref{eq-part-mk} must be negative.
Let $I$ be the smallest index $i$ such that $\mu^{(k)}_{i}< 0$. Note that $M_k$, as a submatrix of $\varphi_t\,\mathrm{JT}_{\lambda}$, has row indices
$\ell(\lambda)-k+1<\cdots<\ell(\lambda).$
Then, for any $\ell(\lambda)-k+I\leq i\leq \ell(\lambda)$, the $i$-th row of $\varphi_t\,\mathrm{JT}_{\lambda}$ must be of the form
$$(0,\ldots,0,1,*,\ldots,*),$$
where a $*$ denotes some nonzero element in $\mathbb{Q}[t]$.
Suppose that the column indices of the $1$'s in these rows are
$j_{I},\ldots,j_k$. Then we must have
$I-1<j_{I}<\cdots<j_k$, i.e., the $1$'s in these rows
appear strictly from left-to-right as we move down $\varphi_t\,\mathrm{JT}_{\lambda}$.
Now consider the $k\times k$ submatrix $N$ of $\varphi_t\,\mathrm{JT}_{\lambda}$ with row indices
$\ell(\lambda)-k+1<\cdots<\ell(\lambda)$ and column indices $1<\cdots<I-1<j_{I}<\cdots<j_k$.
It is easy to check that $N$ is a block matrix of the form
$$
\begin{pmatrix}
\varphi_t\, \mathrm{JT}_{\mu^{(k,+)}} & Q\\
O & P
\end{pmatrix},
$$
where $O$ is a $(k-I+1)\times (I-1)$ zero matrix and $P$ is upper unitriangular. Therefore, $\det N=\varphi_t\,s_{\mu^{(k,+)}}$.

To prove (C2), we only need to consider every $k\times k$ submatrix $N_k$ of $\varphi_t\,\mathrm{JT}_{\lambda}$ with nonzero determinant.
Suppose that $\det N_k=\varphi_t s_{\rho/\sigma}$ for some skew partition $\rho/\sigma$. By Lemma \ref{lem-stanley}, if $\langle s_{\rho/\sigma},\, s_{\tau} \rangle\neq 0$,
then $\mu^{(k,+)}\subseteq \tau$, and hence the contents of $\mu^{(k,+)}$ form a submultiset of the contents of $\tau$. Therefore, $\varphi_t s_{\tau}$ is divisible by $\varphi_t\,s_{\mu^{(k,+)}}$, so
is $\varphi_t s_{\rho/\sigma}=\det N_k$. \qed

We now have all the ingredients necessary for the proof of Theorem \ref{thm-snf-jt}.

\textit{Proof of Theorem \ref{thm-snf-jt}.} From \eqref{eq-decomposition} it follows that
\begin{align}
\varphi_t\,s_{\mu^{(k,+)}}=c_k\prod_{l=1}^k \prod_{(i,j)\in D_{\ell(\lambda)-l+1}}(t+c(i,j)),\label{eq-key-identity}
\end{align}
where $c_k$ is a nonzero rational number. By Proposition \ref{prop-gcd}, the greatest common divisor of the $k\times k$ minors is equal to
$\varphi_t\,s_{\mu^{(k,+)}}$. Combining Theorem \ref{thm-minor} and equation \eqref{eq-key-identity} we obtain the desired result. \qed

\section{The Giambelli matrix} \label{sec-g}

The objective of this section is to give a proof of Theorem \ref{thm-snf-g}. Similar to the proof of Theorem \ref{thm-snf-jt},
we need to compute the greatest common divisor of $k\times k$ minors of $\varphi_t\,\mathrm{G}_{\lambda}$ explicitly. However, the Littlewood-Richardson rule is not necessary for the proof here.

Given a partition $\lambda=(\alpha_1,\ldots,\alpha_r|\beta_1,\ldots,\beta_r)$, it is easy to see that, for $1\leq k \leq r$, the $k$-th diagonal hook of $\lambda$
is  $D_k=(\alpha_k|\beta_k)$.
Let $M_k$ be the square submatrix consisting of the last $k$ rows and last $k$ columns of $\varphi_t\,\mathrm{G}_{\lambda}$.
Let $\mu^{(k)}$ denote the partition $(\alpha_{r-k+1},\ldots,\alpha_r | \beta_{r-k+1},\ldots,\beta_r)$, which has the following hook decomposition
\begin{align}\label{eq-decomposition-g}
\mu^{(k)}=\bigcup_{l=1}^kD_{r-l+1}.
\end{align}
Thus $M_k$ is the specialized Giambelli matrix for the partition $\mu^{(k)}$, and hence
\begin{align}
\det M_k=\varphi_t\,s_{\mu^{(k)}}=c_k\prod_{l=1}^k \prod_{(i,j)\in D_{r-l+1}}(t+c(i,j)),\label{eq-key-identity-g}
\end{align}
where $c_k$ is a nonzero rational number. We have the following result.

\begin{prop}\label{prop-gcd-g}
Let $M_k$ be the square submatrix of the last $k$ rows and last $k$ columns of $\varphi_t\,\mathrm{G}_{\lambda}$. Then every $k\times k$ minor of $\varphi_t\,\mathrm{G}_{\lambda}$ is divisible by $\det M_k$
in the ring $\mathbb{Q}[t]$.
\end{prop}

\proof Consider a $k\times k$ square submatrix
$N_k$ of $\varphi_t\,\mathrm{G}_{\lambda}$ with row indices
$1\leq i_1< \cdots<i_k\leq r$ and column indices
 $1\leq j_1<\cdots<j_k\leq r$.
It is obvious that $N_k$ is the specialized Giambelli matrix for the partition $\nu^{(k)}=(\alpha_{i_1},\ldots,\alpha_{i_k}|\beta_{j_1},\ldots,\beta_{j_k})$.
Since $M_k$ consists of the last $k$ rows and last $k$ columns of $\varphi_t\,\mathrm{G}_{\lambda}$, we have
$\alpha_{i_l}\geq\alpha_{r-k+l}$ and $\beta_{j_l}\geq\beta_{r-k+l}$ for all $1\leq l\leq k$. Thus $\mu^{(k)}\subseteq \nu^{(k)}$,
and hence the contents of $\mu^{(k)}$ form a submultiset of the contents of $\nu^{(k)}$. The desired result immediately follows from the hook-content formula.
\qed

We proceed to finish the proof of Theorem \ref{thm-snf-g}.

\textit{Proof of Theorem \ref{thm-snf-g}.}
By Proposition \ref{prop-gcd-g}, the determinant $\det M_k$ is a greatest common divisor of $k\times k$ minors of $\varphi_t\,\mathrm{G}_{\lambda}$.
The proof then follows from Theorem \ref{thm-minor} and equation \eqref{eq-key-identity-g}.
\qed


\section{The Lascoux-Pragacz matrix}\label{sec-lp}

In this section we will give a proof of Theorem \ref{thm-snf-lp}.
As in the proofs of Theorems \ref{thm-snf-jt} and \ref{thm-snf-g}, we need to
determine the greatest common divisor of $k\times k$ minors of $\varphi_t\,\mathrm{LP}_{\lambda}$.

First, we would like to point out that
each minor of $\varphi_t\,\mathrm{LP}_{\lambda}$ is equal to a skew Schur function (under the specialization $\varphi_t$) by the following result due to Lascoux and Pragacz \cite{lascoux1988ejc}.

\begin{thm}[{\cite[Corallary 4.7]{lascoux1988ejc}}]\label{thm-general-lp}
Let $\lambda$ be a partition of rank $r$, its Frobenius notation  $(\alpha_1,\ldots,\alpha_r|\beta_1,\ldots,\beta_r)$, and its rim decomposition $\theta_1,\ldots,\theta_r$.
Further, let $k$ be a nonnegative integer $\leq r$,
$\mathcal{S}_1=\{i_1,\ldots,i_k\}$, $\mathcal{S}_2=\{j_1,\ldots,j_k\}$ be two subsets of the set $\mathcal{S}=\{1,\ldots,r\}$, and $\mu$ be the partition of Frobenius decomposition $(\alpha_{i_1},\ldots,\alpha_{i_k}|\beta_{j_1},\ldots,\beta_{k_k})$. Then
$$s_{\lambda/\mu}= \det (s_{\theta_i^{+} \& \theta_j^{-}})_{i\in\overline{\mathcal{S}_1},j\in\overline{\mathcal{S}_2}},$$
where $\overline{\mathcal{S}_i}$ is the relative complement of $\mathcal{S}_i$ with respect to $\mathcal{S}$.
\end{thm}

Suppose that $\lambda$ is a partition of rank $r$ with Frobenius decomposition $(\alpha_1,\ldots,\alpha_r|\beta_1,\ldots,\beta_r)$ and rim decomposition $\theta_1,\ldots,\theta_r$.
For $1\leq k \leq r$, let $M_k$ be the square submatrix of the last $k$ rows and last $k$ columns of $\varphi_t\,\mathrm{LP}_{\lambda}$. By Theorem \ref{thm-general-lp},
the determinant $\det M_k$ is equal to $\varphi_t\,s_{\lambda/\mu^{(k)}}$, where $\mu^{(k)}=(\alpha_{1},\ldots,\alpha_{r-k}|\beta_{1},\ldots,\beta_{r-k})$.
It is easy to see that $\lambda/\mu^{(k)}$ is an ordinary partition with Frobenius notation $(\alpha_{r-k+1},
\ldots,\alpha_{r}
|\beta_{r-k+1},\ldots,\beta_{r})$.
Let $N_k$ be the $k\times k$ square submatrix
 of $\varphi_t\,\mathrm{LP}_{\lambda}$ with row indices
$1\leq i_1< \cdots<i_k\leq r$ and column indices
 $1\leq j_1<\cdots<j_k\leq r$. Let $\{i'_1,\ldots,i'_{r-k}\}$ be the relative complement of $\{i_1,\ldots,i_k\}$ with respect to $\{1,\ldots,r\}$,
 and $\{j'_1,\ldots,j'_{r-k}\}$ the relative complement of $\{j_1,\ldots,j_k\}$, ordered increasingly. Again by Theorem \ref{thm-general-lp}, the determinant $\det N_k$ is equal to $\varphi_t\,s_{\lambda/\nu^{(k)}}$, where $\nu^{(k)}=(\alpha_{i'_1},\ldots,\alpha_{i'_{r-k}}|\beta_{j'_1},\ldots,\beta_{j'_{r-k}})$.
 We have the following property, which is similar to Lemma \ref{lem-stanley}.

\begin{lem}\label{lem-lascoux-pragacz} Let $\mu^{(k)}$ and $\nu^{(k)}$ be defined as above.
Then all other cells of
$\lambda/\nu^{(k)}$ are to the left or above of $\lambda/\mu^{(k)}$. Furthermore, we have $\lambda/\mu^{(k)}\subseteq \tau$ if $\langle s_{\lambda/\nu^{(k)}},\, s_{\tau} \rangle> 0$.
\end{lem}

\proof  To prove the first part, it suffices to show that $\nu^{(k)}\subseteq \mu^{(k)}$. This is clear since $\alpha_{i'_l}\leq \alpha_{l}$ and $\beta_{j'_l}\leq \beta_{l}$
for any $1\leq l\leq r-k$.
By the Littlewood-Richardson rule, the statement that $\langle s_{\lambda/\nu^{(k)}},\, s_{\tau} \rangle> 0$ is equivalent to the existence of
a semistandard Young tableau (say $T$) of shape $\lambda/\nu^{(k)}$ and content $\tau$ whose reverse reading word is a lattice permutation.
For simplicity, let $\nu$ denote the ordinary partition $(\lambda/\mu^{(k)})$. We claim that, for each $1\leq i\leq \ell(\nu)$,
the number of occurrences of $i$ in $T$ is at least $\nu_i$. For $i=1$, consider the first row of $\nu$, which lies in the $(r-k+1)$-th row of $T$.
Suppose that this row contains $n_1$ $1$'s, $n_2$ $2$'s, $n_3$ $3$'s, etc. Since $T$ satisfies a lattice permutation condition, the first $r-k+1$ rows of
$T$ must contain at least $(n_1+n_2)$ $1$'s, at least $(n_2+n_3)$ $2$'s, etc. Continuing in this way, there exist at least $(n_1+n_2+n_3+\cdots)$ $1$'s
in the first $r-k+1$ rows of $T$, and hence at least $\nu_1$ $1$'s in $T$. The same reasoning can be used to prove that the tableau $T$ has at least $\nu_2$ $2$'s,
$\nu_3$ $3$'s, etc, keeping in mind that $T$ is increasing along the rows and strictly increasing down the columns.
Therefore, for any $1\leq i\leq \ell(\nu)$, we have $\nu_i\leq \tau_i$, and hence $\nu\subseteq \tau$, as desired.
\qed

Based on the above result, we can determine the greatest common divisor of $k\times k$ minors of $\varphi_t\,\mathrm{LP}_{\lambda}$.

\begin{prop}\label{prop-gcd-lp}
Let $M_k$ be the square submatrix of the last $k$ rows and last $k$ columns of $\varphi_t\,\mathrm{LP}_{\lambda}$.
 Then every $k\times k$ minor of $\varphi_t\,\mathrm{LP}_{\lambda}$ is divisible by $\det M_k$
in the ring $\mathbb{Q}[t]$.
\end{prop}

\proof Let $M_k$ and $N_k$ be chosen as above. By Lemma \ref{lem-lascoux-pragacz}, if $\langle s_{\lambda/\nu^{(k)}},\, s_{\tau} \rangle> 0$,
then $\lambda/\mu^{(k)}\subseteq \tau$, and hence the contents of $\lambda/\mu^{(k)}$ form a submultiset of the contents of $\tau$.
By the hook-content formula $\varphi_t s_{\tau}$ is divisible by $\varphi_t\,s_{\lambda/\mu^{(k)}}=\det M_k$, so
is $\varphi_t s_{\lambda/\nu^{(k)}}=\det N_k$. This completes the proof.
\qed

Now, we are in position to prove Theorem \ref{thm-snf-lp}.

\textit{Proof of Theorem \ref{thm-snf-lp}.}
By Proposition \ref{prop-gcd-lp}, the determinant $\det M_k=\varphi_t\,s_{\lambda/\mu^{(k)}}$ is a greatest common divisor of $k\times k$ minors of $\varphi_t\,\mathrm{LP}_{\lambda}$.
Note that $\lambda/\mu^{(k)}$ has Frobenius decomposition $(\alpha_{r-k+1},
\ldots,\alpha_{r}
|\beta_{r-k+1},\ldots,\beta_{r})$. Thus,
\begin{align*}
\det M_k=\varphi_t\,s_{\lambda/\mu^{(k)}}=c_k\prod_{l=1}^k \prod_{(i,j)\in D_{r-l+1}}(t+c(i,j)).
\end{align*}
for some nonzero rational number $c_k$.
Then the application of  Theorem \ref{thm-minor} establishes the desired result.
\qed

\section{The Giambelli-type matrix}\label{sec-gt}

The main objective of this section is to give a proof of Theorem \ref{thm-snf-gg}.
Let us first review the construction of the Giambelli-type matrices introduced by Hamel and Goulden \cite{hamel1995ejc}.
We will use the reformulation involving the notion of a cutting strip due to Chen, Yan and Yang \cite{chen2005jac}.

Suppose that $\lambda/\mu$ is a connected skew partition with $d$ nonempty diagonals.
Given an outside decomposition $\Pi=(\theta_1,\theta_2,\ldots,\theta_m)$ of $\lambda/\mu$, cells in the same diagonal within the skew diagram
will either all go up or all go right with respect to $\Pi$.
Let $\theta(\Pi)$ be the unique border strip of length $d$, which occupies the same nonempty diagonals as $\lambda/\mu$ and the cells in a given diagonal
goes up or right exactly as the cells of $\lambda/\mu$ do with respect to $\Pi$. The ribbon $\theta(\Pi)$ is called the cutting strip of $\Pi$.
(This establishes a one-to-one correspondence between the set of outside decompositions of $\lambda/\mu$ and the set of border strips with $d$ cells.)
Given any two contents $p,q$, we may define the strip $\theta[p,q]$ as follows:
\begin{itemize}
\item[(1).] If $p\leq q$, then $\theta[p,q]$ is the segment of $\theta(\Pi)$ starting with the cell having content $p$ and ending with the cell having content $q$;
\item[(2).] If $p=q+1$, then $\theta[p,q]$ is the empty strip $\emptyset$ with $s_{\theta[p,q]}=1$;
\item[(3).] If $p>q+1$, then $\theta[p,q]$ is undefined with $s_{\theta[p,q]}=0$.
\end{itemize}
For any strip $\theta_i$ in $\Pi$, let $p(\theta_i)$ denote the content of the starting cell of $\theta_i$, and $q(\theta_i)$ the content of the ending cell of $\theta_i$.
Further let
$$\mathrm{M}_{\lambda/\mu}(\Pi)=(s_{\theta[p(\theta_i),q(\theta_j)]})_{i,j=1}^{m}.$$
This is the Giambelli-type matrix studied by Hamel and Goulden \cite{hamel1995ejc}.
They obtained the following result:
\begin{thm}[{\cite[Theorem 3.1]{hamel1995ejc}}]
For any connected $\lambda/\mu$ and any outside decomposition $\Pi$, we have $s_{\lambda/\mu}=\det \mathrm{M}_{\lambda/\mu}(\Pi)$.
\end{thm}

As we mentioned at the end of the introduction, it might be possible to give a proof of Theorem \ref{thm-snf-gg} by explicitly computing the greatest common divisor of
$k\times k$ minors of $\mathrm{M}_{\lambda}(\Pi)$. 
Here we will prove Theorem \ref{thm-snf-gg} by using the stable equivalence of Giambelli-type matrices, established by Chen and Yang \cite{chen2004ejc}.

The notion of stable equivalence of matrices was introduced by Kuperberg \cite{kuperberg2002ejc}. Let $R$ be a commutative ring with unit.
Let $M$ be an $n\times k$ matrix over
$R$. We say that $M'$ is stably equivalent to $M$ if
$M'$ can be obtained from $M$ under the following operations:
general row operations,
$$M\rightarrowtail AM$$
where $A$ is an $n\times n$ invertible matrix over $R$; general
column operations,
$$M\rightarrowtail MB$$
where $B$ is a $k\times k$ invertible matrix over $R$; and
stabilization
$$M\rightarrowtail
\left(
    \begin{tabular}{c|c}
    1 & 0\\
    \hline
    0 & M
    \end{tabular}
\right)
$$ and its inverse.

Kuperberg \cite[Question 15]{kuperberg2002ejc} asked whether the
Jacobi-Trudi matrix and the dual Jacobi-Trudi matrix are stably equivalent for any skew partition over the ring of symmetric functions.
Chen and Yang \cite{chen2004ejc} proved the following stable equivalence property of Giambelli-type matrices.
\begin{thm}[{\cite[Proposition 3.4]{chen2004ejc}}]\label{thm}
Let $\Pi$ and $\Pi'$ be two outside decompositions of the edgewise connected skew diagram $\lambda/\mu$. Then the Giambelli-type matrices $\mathrm{M}_{\lambda/\mu}(\Pi)$ and $\mathrm{M}_{\lambda/\mu}(\Pi')$
are stably equivalent over the ring of symmetric function.
\end{thm}

Since for any symmetric function $f$ of finite degree the specialization $\varphi_t\,f$ is still a polynomial in $t$,
we immediately have the following result.
\begin{cor}\label{cor-stable-1}
Given any two outside decompositions $\Pi$ and $\Pi'$ of $\lambda/\mu$, the specialized Giambelli-type matrices $\varphi_t\,\mathrm{M}_{\lambda/\mu}(\Pi)$ and $\varphi_t\,\mathrm{M}_{\lambda/\mu}(\Pi')$
are stably equivalent over the polynomial ring $\mathbb{Q}[t]$.
\end{cor}

To prove Theorem \ref{thm-snf-gg}, we also need the following result.

\begin{lem}\label{lem-stable-snf}
Suppose that $M$ is an $n\times n$ nonsingular matrix and $M'$ is an $m\times m$ nonsingular matrix over a principal ideal domain, where
$n\leq m$. If $M$ is stably equivalent to $M'$ and $M\,\stackrel{\mathrm{snf}}{\rightarrow}\,(d_1,d_2,\ldots, d_n)$, then $M'\,\stackrel{\mathrm{snf}}{\rightarrow}\,({1,\ldots,1},d_1,d_2,\ldots, d_n)$, where $1$ occurs $(m-n)$ times before $d_1$.
\end{lem}

\proof The proof is by definition.
\qed

Now we can prove Theorem \ref{thm-snf-gg}.

\textit{Proof of Theorem \ref{thm-snf-gg}.} In Corollary \ref{cor-stable-1}, take $\Pi'$ to be the hook decomposition of $\lambda$, and then $\varphi_t\,\mathrm{M}_{\lambda}(\Pi')$ is the specialized Giambelli matrix $\varphi_t\,\mathrm{G}_{\lambda}$. Combining Lemma \ref{lem-stable-snf} and Theorem \ref{thm-snf-g}, we obtain the desired result, by noting that, for $1\leq k\leq m-r$, the diagonal $D_{m-k+1}$ is the empty partition and hence
\begin{align*}
d_k=\prod_{(i,j)\in D_{m-k+1}}(t+c(i,j))=1.
\end{align*}
\qed

We should mention that Theorem \ref{thm-snf-gg} tells that the Smith normal form of a specialized Giambelli-type matrix can be obtained from that
of the corresponding specialization of the classical Giambelli matrix by a sequence of stabilization operations. Since the Jacobi-Trudi matrix is a special Giambelli-type matrix,
Theorem \ref{thm-snf-jt} can also be considered as a corollary of Theorem \ref{thm-snf-g} and Lemma \ref{lem-stable-snf}. This provides another approach to Stanley's result on the Smith normal form
of a specialized Jacobi-Trudi matrix.

\section{A $q$-analogue} \label{sec-gq}

In \cite{stanley2015arxiv,stanley2016jcta}, Stanley gave a $q$-analogue of Theorem \ref{thm-snf-jt}.
In this section, we shall show that Stanley's approach also enables us to give a $q$-analogue of Theorems \ref{thm-snf-g} and \ref{thm-snf-lp}.

Let us first give an overview of Stanley's $q$-analogue of Theorem \ref{thm-snf-jt}.
Letting $t$ be a fixed positive integer, a natural $q$-analogue $\varphi_t(q)\,f$ of $\varphi_t\,f$ is given by
$f(1,q,\ldots,q^{t-1})$. It is well known that
\begin{align}
\varphi_t(q)\,s_{\lambda}={q^{\sum(k-1)\lambda_k}}\prod_{(i,j)\in \lambda}\frac{1-q^{t+c(i,j)}}{1-q^{h(i,j)}}.
\end{align}
Then setting $y=q^t$, the specialization $\varphi_t(q)\,f$ for any symmetric function $f$ becomes a polynomial in $y$ with coefficients in the field $\mathbb{Q}(q)$.
Let $\widehat{\mathrm{JT}}_{\lambda}$ denote the matrix obtained from $\varphi_t(q)\,{\mathrm{JT}}_{\lambda}$ by substituting $q^t=y$.
Stanley \cite{stanley2015arxiv} studied the Smith normal form of $\widehat{\mathrm{JT}}_{\lambda}$ over the principal ideal domain $\mathbb{Q}(q)[y]$, and obtained the following result.

\begin{thm}[{\cite[Theorem 3.1]{stanley2015arxiv}}]\label{thm-snf-jt-q1}
Suppose that $\widehat{\mathrm{JT}}_{\lambda}\,\stackrel{\mathrm{snf}}{\rightarrow}\, (d_1,d_2,\ldots,d_{\ell(\lambda)})$ over $\mathbb{Q}(q)[y]$.
Then, for $1\leq k\leq \ell(\lambda)$, we can take
\begin{align*}
d_k=\prod_{(i,j)\in D_{\ell(\lambda)-k+1}}(1-q^{c(i,j)}y).
\end{align*}
\end{thm}

Stanley remarked that the proof of Theorem \ref{thm-snf-jt} given in Section \ref{sec-jt} carries over to Theorem \ref{thm-snf-jt-q1}.
In the same manner, we can give a $q$-analogue of Theorem \ref{thm-snf-g}. The proof of Theorem \ref{thm-snf-g} given in Section \ref{sec-g} carries over to this $q$-version.

\begin{thm}\label{thm-snf-g-q1}
Let $\widehat{\mathrm{G}}_{\lambda}$ denote the matrix obtained from $\varphi_t(q)\,{\mathrm{G}}_{\lambda}$ by substituting $q^t=y$.
Suppose that $\widehat{\mathrm{G}}_{\lambda}\,\stackrel{\mathrm{snf}}{\rightarrow}\, (d_1,d_2,\ldots,d_{r})$ over $\mathbb{Q}(q)[y]$, where $r=\mathrm{rank}(\lambda)$.
Then, for $1\leq k\leq r$, we can take
\begin{align*}
d_k=\prod_{(i,j)\in D_{r-k+1}}(1-q^{c(i,j)}y).
\end{align*}
\end{thm}

We can also give a $q$-analogue of Theorem \ref{thm-snf-lp}.

\begin{thm}\label{thm-snf-lp-q1}
Let $\widehat{\mathrm{LP}}_{\lambda}$ denote the matrix obtained from $\varphi_t(q)\,{\mathrm{LP}}_{\lambda}$ by substituting $q^t=y$.
Suppose that $\widehat{\mathrm{LP}}_{\lambda}\,\stackrel{\mathrm{snf}}{\rightarrow}\, (d_1,d_2,\ldots,d_{r})$ over $\mathbb{Q}(q)[y]$, where $r=\mathrm{rank}(\lambda)$.
Then, for $1\leq k\leq r$, we can take
\begin{align*}
d_k=\prod_{(i,j)\in D_{r-k+1}}(1-q^{c(i,j)}y).
\end{align*}
\end{thm}

Next we will give a $q$-analogue of Theorem \ref{thm-snf-gg}. Given an outside decomposition $\Pi$ of $\lambda$, let $\widehat{\mathrm{M}}_{\lambda}(\Pi)$
denote the matrix obtained from $\varphi_t(q)\,{\mathrm{M}}_{\lambda}(\Pi)$ by substituting $q^t=y$. We have the following result.

\begin{thm}\label{thm-snf-gg-q1}
Suppose that ${\mathrm{M}}_{\lambda}(\Pi)$ is of order $m$ and $\widehat{\mathrm{M}}_{\lambda}(\Pi)\,\stackrel{\mathrm{snf}}{\rightarrow}\, (d_1,d_2,\ldots,d_{m})$ over $\mathbb{Q}(q)[y]$.
Then, for $1\leq k\leq m$, we can take
\begin{align*}
d_k=\prod_{(i,j)\in D_{m-k+1}}(1-q^{c(i,j)}y).
\end{align*}
\end{thm}

\proof The proof of Theorem \ref{thm-snf-gg} given in Section \ref{sec-gt} carries over to this $q$-analogue. Given any two outside decompositions $\Pi$ and $\Pi'$ of $\lambda/\mu$, Theorem \ref{thm} tells that the matrices $\mathrm{M}_{\lambda/\mu}(\Pi)$ and $\mathrm{M}_{\lambda/\mu}(\Pi')$
are stably equivalent over the ring of symmetric functions. This implies that $\widehat{\mathrm{M}}_{\lambda/\mu}(\Pi)$ and $\widehat{\mathrm{M}}_{\lambda/\mu}(\Pi')$
are stably equivalent over the polynomial ring $\mathbb{Q}(q)[y]$. In particular, we obtain the stable equivalence of $\widehat{\mathrm{M}}_{\lambda}(\Pi)$ and $\widehat{\mathrm{G}}_{\lambda}$.
The desired result immediately follows from Lemma \ref{lem-stable-snf}. \qed

However, the above $q$-analogue seems unsatisfactory since one can not reduce to the original matrices by substituting $y=1$. Stanley \cite{stanley2015arxiv} suggested another $q$-analogue of
Theorem \ref{thm-snf-jt} in the following way. For any symmetric function $f$, let $\varphi^{*}\,f$ denote the substitution $q^t\rightarrow \frac{1}{(1-q)y+1}$ after writing $\varphi_t(q)\,f$ as a
polynomial in $q$ and $q^t$. Under this substitution, for any $k\in\mathbb{Z}$ there holds
$$
1-q^ky\rightarrow \frac{(1-q)(y+\mathbf{(k)})}{(1-q)y+1},
$$
where $\mathbf{(k)}=\frac{1-q^k}{1-q}$ for any $k\in\mathbb{Z}$. Thus, for any partition $\lambda$, we have
\begin{align*}
\varphi^* s_{\lambda}
&=q^{\sum (k-1)\lambda_k}\prod_{(i,j)\in \lambda}\frac{1-q^{c(i,j)}\cdot\frac{1}{(1-q)y+1}}{1-q^{h(i,j)}}\\
&=\frac{q^{\sum (k-1)\lambda_k}}{((1-q)y+1)^{|\lambda|}}\prod_{(i,j)\in \lambda} \frac{(y+\textbf{c(i,j)})}{\textbf{(h(i,j)})}.
\end{align*}
Then, we define a map $\varphi^{\diamond}$ from the ring of symmetric functions $\Lambda$ to $\mathbb{Q}(q)[y]$ by letting
\begin{align}
\varphi^{\diamond} s_{\lambda}&= ((1-q)y+1)^{|\lambda|}\varphi^* s_{\lambda}=q^{b(\lambda)}\prod_{(i,j)\in \lambda} \frac{(y+\textbf{c(i,j)})}{\textbf{(h(i,j)})}\label{q-analogue}
\end{align}
and then extend linearly. It is easy to show that $\varphi^{\diamond} (s_{\lambda}\cdot s_{\mu})=\varphi^{\diamond} s_{\lambda}\cdot\varphi^{\diamond} s_{\mu}$.
Therefore, $\varphi^{\diamond}$ is well-defined, and it is a homomorphism from $\Lambda$ to $\mathbb{Q}(q)[y]$.
Then $\widetilde{\mathrm{JT}}_{\lambda}=\varphi^{\diamond}\,{\mathrm{JT}}_{\lambda}$ is just the specialized Jacobi-Trudi matrix studied by Stanley \cite{stanley2015arxiv},
who obtained the following result.

\begin{thm}[{\cite[Theorem 3.2]{stanley2015arxiv}}]\label{thm-snf-jt-q2}
Suppose that $\widetilde{\mathrm{JT}}_{\lambda}\,\stackrel{\mathrm{snf}}{\rightarrow}\, (d_1,d_2,\ldots,d_{\ell(\lambda)})$ over $\mathbb{Q}(q)[y]$.
Then, for $1\leq k\leq \ell(\lambda)$, we can take
\begin{align*}
d_k=\prod_{(i,j)\in D_{\ell(\lambda)-k+1}}(y+\textbf{c(i,j)}).
\end{align*}
\end{thm}

For the Giambelli matrix and the Lascoux-Pragacz matrix, let $\widetilde{\mathrm{G}}_{\lambda}=\varphi^{\diamond}\,{\mathrm{G}}_{\lambda}$ and $\widetilde{\mathrm{LP}}_{\lambda}=\varphi^{\diamond}\,{\mathrm{LP}}_{\lambda}$.
Similar arguments enable us to get a natural $q$-analogue of Theorem \ref{thm-snf-g} and Theorem \ref{thm-snf-lp}.

\begin{thm}\label{thm-snf-g-q2}
Suppose that $\widetilde{\mathrm{G}}_{\lambda}\,\stackrel{\mathrm{snf}}{\rightarrow}\, (d_1,d_2,\ldots,d_{r})$ over $\mathbb{Q}(q)[y]$, where $r=\mathrm{rank}(\lambda)$.
Then, for $1\leq k\leq r$, we can take
\begin{align*}
d_k=\prod_{(i,j)\in D_{r-k+1}}(1-q^{c(i,j)}y).
\end{align*}
\end{thm}

\begin{thm}\label{thm-snf-lp-q2}
Suppose that $\widetilde{\mathrm{LP}}_{\lambda}\,\stackrel{\mathrm{snf}}{\rightarrow}\, (d_1,d_2,\ldots,d_{r})$ over $\mathbb{Q}(q)[y]$, where $r=\mathrm{rank}(\lambda)$.
Then, for $1\leq k\leq r$, we can take
\begin{align*}
d_k=\prod_{(i,j)\in D_{r-k+1}}(1-q^{c(i,j)}y).
\end{align*}
\end{thm}

Theorem \ref{thm-snf-gg} also enjoys a $q$-analogue under the specialization $\varphi^{\diamond}$.
Given an outside decomposition $\Pi$ of $\lambda$, let $\widetilde{\mathrm{M}}_{\lambda}(\Pi)=\varphi^{\diamond}\,{\mathrm{M}}_{\lambda}(\Pi)$.
Noting that $\varphi^{\diamond}$ is a homomorphism from $\Lambda$ to $\mathbb{Q}(q)[y]$, the proof of Theorem \ref
{thm-snf-gg-q1} carries over directly to the following result.

\begin{thm}\label{thm-snf-gg-q2}
Suppose that ${\mathrm{M}}_{\lambda}(\Pi)$ is of order $m$ and $\widetilde{\mathrm{M}}_{\lambda}(\Pi)\,\stackrel{\mathrm{snf}}{\rightarrow}\, (d_1,d_2,\ldots,d_{m})$ over $\mathbb{Q}(q)[y]$.
Then, for $1\leq k\leq m$, we can take
\begin{align*}
d_k=\prod_{(i,j)\in D_{m-k+1}}(y+\textbf{c(i,j)}).
\end{align*}
\end{thm}

\noindent{\bf Acknowledgements.} This work was supported by the National Science Foundation of China.


\begin{thebibliography}{1}


\bibitem{chen2004ejc}
W. Y. C. Chen  and A. L. B. Yang, Stable equivalence over symmetric functions,
\emph{Electron. J. Combin}. \textbf{11}(2) (2004--2006), \#R23.


\bibitem{chen2005jac}
W. Y. C. Chen, G. G. Yan  and A. L. B. Yang, Transformations of border strips and schur function determinants, \emph{J. Algebraic Combin}. \textbf{21} (2005), 379--394.

\bibitem{giambelli1903at} G. Z. Giambelli, Alcune propriet\`a delle funzioni simmetriche caratteristiche,
Atti Torino. \textbf{38} (1903), 823--844.


\bibitem{hamel1995ejc}
A. M. Hamel  and I. P. Goulden, Planar decompositions of tablaux and Schur function determinants, \emph{Europ. J. Combin.} \textbf{16} (1995), 461--477.

\bibitem{kuperberg2002ejc}
G. Kuperberg, Kasteleyn cokernels, \emph{Electron. J. Combin.} \textbf{9}(1) (2002), \#R29.

\bibitem{lascoux1988ejc}
A. Lascoux  and P. Pragacz, Ribbon Schur function, \emph{Europ. J. Combin.} \textbf{9} (1988), 561--574.

\bibitem{stanley1999ec2}
R. P. Stanley, {Enumerative Combinatorics}, Vol. 2, Cambridge University Press, Cambridge, 1999.

\bibitem{stanley2015arxiv}
R. P. Stanley, The Smith normal form of a specialized Jacobi-Trudi matrix, \emph{Europ. J. Combin.} \textbf{62} (2017), 178--182.

\bibitem{stanley2016jcta}
R. P. Stanley, Smith normal form in combinatorics, \emph{J. Combin. Theory, Ser. A} \textbf{144} (2016), 476--495.


%
%
%
%
%
%
%
%
%
%


\end{thebibliography}
\end{document}